\documentclass[12pt]{article}
\usepackage[english]{babel}
\usepackage{amssymb,amsmath,graphics,tikz}

\usepackage{hyperref}

\parskip\medskipamount
\newtheorem{theorem}{Theorem}[section]

\newtheorem{prop}[theorem]{Proposition}

\def\<{\langle}
\def\>{\rangle}
\newcommand{\proof}{\emph{Proof.~}}

\def\qed{{\hfill\hphantom{.}\nobreak\hfill$\Box$}}

\newcommand{\define}{\mathrel{\mathop:}=}

\newcommand{\cC}{\mathcal{C}} 

\begin{document}

\author{Koen Struyve}
\title{Free constructions of geometries of Coxeter type}

\maketitle
\begin{abstract}
We establish two free constructions of geometries of Coxeter type. The first construction deals with any Coxeter diagram having no subdiagram of type $\mathsf{A}_3$, the second one with diagrams of type $\mathsf{C}_n$ and $\mathsf{H}_4$.
\end{abstract}

\section{Introduction}
In~\cite{Tit:81} Jacques Tits develops a local approach to buildings, in terms of geometries of type $M$, where $M$ is a Coxeter diagram. 
These are motivated by the applications to the study of sporadic groups, see for example~\cite{Tit:80}.  
While buildings of type $M$ form a subclass of these geometries of type $M$, the converse is not true, as exhibited by the Neumaier geometry (see~\cite{Neu:84}).

One of the results in~\cite{Tit:81} shows that geometries of type $\mathsf{A}_n$ correspond exactly with $n$-dimensional projective geometries, and that every geometry of type $M$ is 2-covered by a building provided that its residues of type $\mathsf{C}_3$ and $\mathsf{H}_3$ are. 
In particular the latter is true in the absence of such subdiagrams.

Tits also mentions that it is possible to construct geometries `freely' whenever $M$ has no subdiagram of type $\mathsf{A}_3$.
For this he refers to a forthcoming paper, which however remained unpublished. 
One of the main motivations for such a construction would be that it, by passing to the universal cover, yields free constructions of buildings for diagrams $M$ with no subdiagram of type $\mathsf{A}_3$, $\mathsf{C}_3$, or $\mathsf{H}_3$.
A more elegant and direct solution to this problem was obtained somewhat later by Mark Ronan (see~\cite{Ron:86}).

The obstruction concerning subdiagrams of type $\mathsf{A}_3$ is a natural one, as these correspond to 3-dimensional projective spaces.

In the last few years there has been renewed interest in geometries of Coxeter type, as they appear naturally in the theory of polar actions on simply connected positively curved manifolds, see~\cite{Fan-Gro-Tho:*}. In particular polar actions on the Cayley plane give rise to two previously unknown geometries of type $\mathsf{C}_3$ not covered by a building. (See also~\cite{Kra-Lyt:14}.)

In this paper we set out to provide a free construction for geometries of type $M$, whenever $M$ has no subdiagram of type $\mathsf{A}_3$, as alluded to by Tits. We also provide constructions of geometries of type $\mathsf{C}_n$ and $\mathsf{H}_n$ where the residues of type $\mathsf{A}_{n-1}$ are all isomorphic to a chosen countably infinite building of this type. 
These results provide a positive answer to Problem 1 of~\cite{Kra-Lyt:14}. Additionally we show how to construct highly homogeneous geometries of spherical type.

The paper is organized as follows. In Section~\ref{section:def} we review some basic notions concerning geometries. In Section~\ref{section:A3} we provide a free construction provided there is no subdiagram of type $\mathsf{A}_3$, and in Section~\ref{section:Cn} we discuss a construction for a class of diagrams including the diagrams $\mathsf{C}_n$ and $\mathsf{H}_4$. Finally, in Section~\ref{section:hom}, Fra\"iss\'e limits are used to construct homogeneous geometries.

\textbf{Acknowledgement.  -- } The author wants to thank Alexander Lytchak for his hospitality and suggesting the problem considered in this paper.
\section{Definitions}\label{section:def}
In this section we introduce some definitions which we will use later on.
\subsection{Geometries}
A \emph{geometry over $I$} is a system $\Gamma \define (V,\tau, *)$, consisting of a set $V$, a map $\tau: V \to I$, and a binary symmetric relation $*$ on $V$ such that for any two elements $x,y$ of $V$ whose images under $\tau$ are identical, the relation $x * y $ holds if and only if $x=y$. The relation `$*$' is the \emph{incidence relation}, the image by $\tau$ of an element or a subset of $V$ is its \emph{type}.

A \emph{flag} of the geometry $\Gamma$ is a set of pairwise incident elements of $V$. Two flags, or a flag and an element of $V$, are said to be incident if their union is a flag. The \emph{rank} of a flag $X$ is its cardinality, the \emph{corank} the cardinality of $I \setminus \tau(X)$.

Let $X$ be a flag, and let $Y$ be the set of all elements of $V \setminus X$ incident to $X$. 
The system $\Gamma_X := (Y, \tau\vert_Y, * \cap (Y \times Y))$ forms a geometry over $I \setminus \tau(X)$ and is called the \emph{residue of $X$} in $\Gamma$.

The geometry $\Gamma$ is \emph{connected} if the graph with vertices $V$ and adjacency relation `$*$' is connected. It is \emph{residually connected} if the residue of every flag of corank 2 is connected.

A geometry is \emph{thick} if every flag of corank 1 is contained in at least three maximal flags. 

\subsection{Rank 2 geometries}
Geometries of rank 2 are basically bipartite graphs. We say that a rank 2 geometry is \emph{without $k$-gons} if the associated bipartite graph does not contain any cycles of length $2k$.

A \emph{generalized $n$-gon} (with $2 \leq n < \infty$) is a thick rank 2 geometry such that the associated graph has girth $2n$ (i.e. the smallest cycle has length $2n$) and diameter $n$. A generalized $\infty$-gon is a thick connected rank 2 geometry without cycles.

\subsection{Geometries of type $M$}
Let $M$ be a Coxeter diagram over $I$ (with associated Coxeter matrix $(m_{i,j})_{i,j \in I}$). We define the \emph{geometries of type $M$} as the thick residually connected geometries over $I$, such that for every $i,j \in I$ ($i \neq j$), the residue of any flag of type $I \setminus \{i,j\}$ is a generalized $m_{i,j}$-gon.


\section{Geometries without $\mathsf{A}_3$ subdiagrams}\label{section:A3}

Fix a connected Coxeter diagram $M$ over $I$ without subdiagrams of type $\mathsf{A}_3$. The goal of this section is to construct a geometry of type $M$ via free construction.

To do so we make use of geometries satisfying a set of specific properties outlined in Section~\ref{section:inter}. We then proceed by showing that one can apply certain extension procedures (see Section~\ref{section:ext}) to such geometries in such a way that these properties are preserved. 
In Section~\ref{section:constr} we obtain the desired geometry of type $M$ by considering the direct limit of a sequence of geometries obtained by these procedures. 

\subsection{Intermediate objects}\label{section:inter}
As intermediate objects in our construction we consider geometries $\Gamma \define (V,\tau,*)$ over $I$ (so $\tau$ maps $V$ to $I$), satisfying the following three properties.
\begin{itemize}
\item[(F)] If $v_1$ and $v_2$ are two vertices such that $\tau(v_1)$ and $\tau(v_2)$ are not adjacent in the Coxeter diagram, then we always have that $v_1 * v_2$.
\item[(P)]
Let $i,j \in I$ be two adjacent vertices of the Coxeter diagram $M$, and $X$ a flag of type $J$ where $i,j \notin J$ and such that every vertex in $I \setminus \{i,j\}$ adjacent to $i$ or $j$ is contained in $J$, then the restriction of the residue of this flag to the vertices of types $i$ and $j$ is a rank 2 geometry without $t$-gons for $t < m_{i,j}$.
\item[(D)] 
If $i,j \in I$ ($i \neq j$) with $m_{i,j} \geq 4$, then the restriction of the geometry $\Gamma$ to the vertices of types $i$ and $j$ is without digons. 
\end{itemize}
Note that the class of such objects is closed under direct limits. The (F) stands for `flat', (P) for `partial', (D) for `digons'.

\subsection{Extension procedures}\label{section:ext}
In this section we describe three procedures to extend a geometry $\Gamma \define (V,\tau,*)$ over $I$ satisfying Properties (F), (P) and (D) to a new geometry $\Gamma'$  over $I$ still satisfying these properties.

\paragraph{Procedure A -- Completing flags.}
Let $X$ be some flag in $\Gamma$ (which we allow to be empty). Let $i$ be a type in $ I\setminus \tau(X)$.
We then add a single vertex $x$ of type $i$ to our geometry.

We also extend the incidence relation $*$ with incidences of the vertex $x$ to certain vertices in $V$, specifically to
\begin{enumerate}
\item
the vertices of the flag $X$, 
\item
and to those vertices for which the incidence would be guaranteed by Property (F).
\end{enumerate}

We obtain a new geometry $(V',\tau,*)$ by this procedure (with $V' = V \cup \{x\}$ and slight abuse of notation for $\tau$ and `$*$'), where the flag $X$ is contained in the larger flag $X \cup \{x\}$.

The geometry $(V',\tau,*)$ satisfies Property (F) by construction. 
In order to verify Properties (P) and (D) notice that for any type $j$ adjacent to $i$, the vertex $x$ is adjacent with at most one vertex of type $j$. 

as the point $x$ is not incident with two vertices of a common type $j$ adjacent to $i$, so $x$ cannot give rise to digons or other $k$-gons.

\paragraph{Procedure B -- Adding paths.}
Let $X$, $J$, $i$ and $j$, be as in the statement of Property (P) with $m:= m_{i,j} < \infty$. 

Let $x$ and $y$ be two vertices of types $i$ or $j$, being of the same type if $m$ is odd, of different type if $m$ is even, and at distance at least $m+1$ in the residue of $X$ restricted to vertices of types $i$ and $j$. This distance is allowed to be infinite.

We extend the vertex set $V$ by adding a path $\gamma \define (x_1 \define x, x_2, \dots, x_{m-1}, x_{m} \define y)$ of length $m-1$ between both, where the types alternate between $i$ and $j$. The incidence relation is extended such that the  added vertices are only incident to
\begin{enumerate}
\item
their adjacent vertices in the path $\gamma$, 
\item 
the vertices of the flag $X$, 
\item
and those vertices for which the incidence would be guaranteed by Property (F).
\end{enumerate}

We denote the newly obtained geometry by $(V',\tau, *)$ (again with abuse of notation), and the subpath of $\gamma$ consisting of the newly added vertices by $\gamma'$.

We postpone the verification of Properties (F), (P) and (D) for the new geometry till Proposition~\ref{prop:ver}.

\paragraph{Procedure C -- Connecting residues.}
Let $X$ be some flag of our geometry of corank at least two with a disconnected residue. 

Let $i,j $ be two different types in $I \setminus \tau(X)$. 

Let $x$ and $y$ be two vertices of type $i$ or $j$ in different connected components of the residue of $X$. 
We now extend the geometry by a new path $\gamma$ where the types alternate between $i$ and $j$, starting at $x$, ending at $y$ and of length at least 4.

The newly added vertices are incident to
\begin{enumerate}
\item
their adjacent vertices in the path $\gamma$, 
\item 
the vertices of the flag $X$, 
\item
and those vertices for which the incidence would be guaranteed by Property (F).
\end{enumerate}

The newly obtained geometry is again denoted by $(V',\tau, *)$. We now verify Properties (F), (P) and (D) for  the last two procedures.

\begin{prop}\label{prop:ver}
The geometries $(V',\tau,*)$ obtained by Procedures B and C still satisfy Properties (F), (P) and (D).
\end{prop}
\proof
Property (F) is satisfied directly by construction.

Next we verify Property (D). 
Suppose, by way of contradition, that our additions give rise to a digon with vertices of types $k$ and $l$ (where $m_{k,l} \geq 4$). 
Without loss of generality we can assume that $k=i$.
Let $z$ be a vertex of type $i$ in the digon. 
As $j$ is the only type adjacent to $i$ for which there are two vertices incident with $x$, we have that $l = j$.

However, if one solely considers incidences between vertices of types $i$ and $j$, then the only thing we did was adding a path of length at least 3 (as $m_{k,l} \geq  4$) between two existing vertices of these types. 
Therefore it is impossible to have created a digon and we conclude that $(V',\tau,*)$ satisfies Property (D).

The last property to verify is Property (P). 
Let $Y$, $k$ and $l$ be the flag and vertices of the Coxeter diagram as required in the statement of Property (P) (so $m_{k,l}\geq 3$). 
We now prove that Property (P) is satisfied for this flag and these types. 
We may assume without loss of generality that $Y$ is minimal (so the types occurring in $Y$ are exactly those adjacent to $k$ or $l$).

First consider the case where $Y$ contains a newly added vertex, for example a vertex $z$ of type $i$. 
The only types for which there exist at least two vertices incident to $z$ are the types $j$, and the types not adjacent to $i$.
So we can assume without loss of generality that $k = j$ and that  $l$ is a type not adjacent to $i$. 
As $z$ is incident with exactly two vertices of type $j$, any $t$-gon on the types $i$ and $j$ in the residue of $Y$ is hence a digon, implying that $m_{k,l} = 3$ by Property (D).
Also $m_{i,j} =3$, as otherwise the two vertices of type $j$ incident with $z$ cannot be incident with more than one vertex of type $l$.
However, this would imply that the types $i$, $j$ and $l$ form a subdiagram of type $\mathsf{A}_3$, which contradicts our assumption on $M$.

Secondly assume that we are in the case that the flag $Y$ is completely contained in $V$. 
If a $t$-gon formed by the vertices of types $k$ and $l$ would contain a newly-added vertex $z$ of type $i$ (so w.l.o.g. $i=k$), then this $t$-gon also contains a vertex of type $j$. 
This as the only type adjacent to $i$ for which $z$ has two incident vertices $a$ and $b$ is the type $j$. 
So we know that $j = l$.

We now claim that every vertex of the flag $X$ for which the type is adjacent to the types $i$ or $j$ is also contained in the flag $Y$. Assume by way of contradiction that this is not the case for a vertex $u$ of the flag $X$ of type $h$, where $h$ is adjacent to either $i$ or $j$.
Denote by $v$ the vertex of type $h$ in the flag $Y$.
As the vertex $z$ is incident with both $u$ and $v$, one has that the types $i$ and $h$ are not adjacent. 
If one of the $a$ or $b$ would also be newly-added (so in $V' \setminus V$), then the same reasoning would also yield that $j$ and $h$ are not adjacent, which is a contradiction. 
So both $a$ and $b$ are already in $V$, implying that $m_{i,j}=3$. The vertices $a,b,u$ and $v$ now form a digon, implying that also $m_{j,h}=3$, which again leads to a contradiction due to the non-existence of subdiagramss of type $\mathsf{A}_3$.

As we now know that every vertex of the flag $X$ for which the type is adjacent to $i$ or $j$ is also a vertex of the flag $Y$, Property (P) for the flag $Y$ and types $k$, $l$ follows from construction.
\qed

\subsection{Construction of a geometry of type $M$}\label{section:constr}

In order to construct a geometry of type $M$ we start with a geometry $\Delta_0 \define (V,\tau,*)$ satisfying the properties outlined in Section~\ref{section:inter}.  
In particular this goemetry may be empty

By applying the Procedures A, B and C outlined in Section~\ref{section:ext} to every possible (viable) combination of elements in $\Delta_0$ and taking the direct limit we obtain a geometry $\Delta_1$. Repeating this step yields a sequence of geometries $\Delta_0, \Delta_1, \Delta_2, \dots$ for which the direct limit $\Delta_\omega$ has the following properties:

\begin{itemize}
\item Every flag is contained a maximal flag of size $\vert I \vert$.
\item The residue of a flag of type $I \setminus \{i,j\}$ (with $i \neq j$), is a generalized $n$-gon with $n \define m_{i,j}$. (Notice that thickness is implied by the repeated application of Procedure A to non-maximal flags.)
\item Every flag of corank at least two has a connected residue.
\end{itemize}

We hence have constructed a geometry of type $M$, exactly what we had set out to do.

\section{Geometries of type $\mathsf{C}_n$ and $\mathsf{H}_4$}\label{section:Cn}

The previous construction does not apply to the types $\mathsf{C}_n$ ($n \geq 4$) and $\mathsf{H}_4$ as these contain subdiagrams of type $\mathsf{A}_3$. The heuristic reason for the exclusion of such subdiagrams is that geometries of type $\mathsf{A}_3$ are exactly the three-dimensional projective spaces defined over (skew) fields, which cannot result from free construction.

In this section we circumvent this restriction by starting with a given projective space and building a geometry of the desired type from it.

\subsection{Setting}\label{section:setting}
Let $M$ be the following Coxeter diagram on $n$ nodes where $n \geq 3$ and $\infty > m \geq 4$.

\begin{center} \begin{tikzpicture}[style=thick]
\foreach \x in {0,1}{
\fill (\x, 0) circle (2pt);
}
\draw (0,0) -- (1.4,0);
\draw[style=dotted] (1.5,0) -- (2.0,0);
\fill (2.5, 0) circle (2pt);
\fill (3.5, 0) circle (2pt);
\draw (2.1,0) -- (3.5,0);

\draw (3,.25) node {$m$} ;
\end{tikzpicture}
\end{center}

The diagrams $\mathsf{C}_n$ ($n \geq 3$) and $\mathsf{H}_n$ ($n= 3$ or $4$) are examples of such diagrams.

Let $\Gamma$ be an vector space of dimension $n$ defined over a (skew) field of (infinite) countable cardinality.

Fix $\Delta_0 \define (V, \tau, *)$ to be the geometry with type set $I \define \{1,\dots, n\}$, where the vertices of type $i \in I$, with $i < n$, are the sub-vector spaces of $\Gamma$ of dimension $i$. 
There are no vertices of type $n$ (yet). 
Two vertices are incident if the corresponding sub-spaces are nested.

\subsection{Intermediate objects}\label{section:inter2}
The intermediate objects in our construction will be geometries $\Delta :=  (W, \tau, *)$ with type set $I$, containing $\Delta_0$, and satisfying the following six properties.
\begin{itemize}
\item[(F)]
Two vertices of types $i$ and $n$, where $i \leq n-2$, are always incident.
\item[(I)]
The number of vertices of type $n$ is finite, the set of all vertices is countable.
\item[(V)]
The type of a vertex in $W \setminus V$ is either $n-1$ or $n$.
\item[(P)]
If $x$ is a vertex of type $n-2$, then the vertices of types $n-1$ and $n$ incident with $y$ form a rank 2 geometry without $k$-gons for $k<m$.
\item[(H)]
If $x$ is a vertex of type $n-1$, then the set of vertices of types $\{1,\dots, n-2\}$ incident with $x$ matches the set of vertices of these types incident with a certain vertex $y$ of type $n-1$ in $\Delta_0$. We say that $y$ is the \emph{precursor} of the vertex $x$.
\item[(C)]
If $x$ is a vertex of type $n$ in $\Delta$, and $y$ is a vertex of type $n-1$ in $\Delta_0$, then there is a unique vertex of type $n-1$ with precursor $y$ incident with $x$.
\end{itemize}
Note that the combination of Properties (F) and (C) implies that the residue of each vertex of type $n$ is isomorphic to the building of type $\mathsf{A}_{n-1}$ associated to the projective space defined on $\Gamma$.

\subsection{Extension procedure}\label{section:ext2}
Let $\Delta :=  (W, \tau, *)$ be a geometry satisfying the properties listed in Section~\ref{section:inter2}.

Let a vertex $z$ of type $n-2$. 
By Property (P) the vertices of types $n-1$ and $n$ incident with $z$ form a rank 2 geometry without $k$-gons for $k < m$. 
Let $x$ and $y$ be two such vertices of types $n-1$ or $n$ incident with $z$, at distance at least $m+1$ (which may be infinite) from each other in this rank 2 geometry, and such that their types agree if $m$ is odd and are different if $m$ is even.

If $x$ is of type $n$, then we pick a vertex $x_1$ of type $n-1$ incident with both $x$ and $z$ such that $x$ is the unique vertex of type $n$ incident with $x_1$.
We claim that is always possible.
Set $K$ to be set of vertices of type $n-1$ incident with both $x$ and $z$. 
Observe that $K$ is infinite by Property (C), and that every vertex of type $n$ different from $x$ is incident with at most one vertex in $K$, as otherwise there would be digons contradicting Property (P). 
The claim then follows by the fact that there are only a finite number of vertices of type $n$ by Property (I).

If $y$ is of type $n$, we pick a vertex $x_{m-2}$ incident with both $y$ and $z$ in an analogous way.

We now extend the vertex set $V$ by adding a path $\gamma \define (x_0 \define x, x_1, \dots, x_{m-2}, x_{m-1} \define y)$ of length $m-1$, where the types alternate between $n-1$ and $n$. (If applicable, the vertices $x_1$ and $x_{m-1}$ are as previously defined.)

For each of the newly added vertices $x_i$ of type $n-1$ we pick an arbitrary precursor incident with $z$, in such a way that the precursors are pairwise different. 
Such a newly added vertex of type $n-1$ is incident only to

\begin{enumerate}
\item
its adjacent vertices (of type $n$) in the path $\gamma$,
\item
the vertices of type $j \leq n-2$ which are incident with the precursor of the newly added vertex.
\end{enumerate}

We also add the following vertices of type $n-1$. 
For each newly added vertex $x_i$ of type $n$ and vertex $a$ of type $n-1$ in $\Delta_0$, which is not a precursor of either $x_{i-1}$ or $x_{i+1}$, we add a vertex of type $n-1$ with precursor $a$ to the geometry. 
This vertex is then incident to 
\begin{enumerate}\setcounter{enumi}{2}
\item
the vertex $x_i$ of type $n$,
\item
the vertices of type $j \leq n-2$ which are incident with the precursor of the newly added vertex.
\end{enumerate}

We end by listing the incidences for a newly added vertex of type $n$.
\begin{enumerate}\setcounter{enumi}{4}
\item
their adjacent vertices (of type $n-1$) in the path $\gamma$,
\item
every vertex of type $j \leq n - 2$,
\item
the incidences already covered in item 3.
\end{enumerate}

We denote the newly obtained geometry by $(W', \tau, *)$. 
The next proposition validates the extension procedure.

\begin{prop}
The geometry $(W', \tau, *)$ obtained by the extension procedure satisfies the properties listed in Section~\ref{section:inter2}.
\end{prop}
\proof
Property (F) is satisfied by construction (see item 6 on the list of incidences).

As we only add a finite number of elements of type $n$, Property (I) is automatically satisfied. 
Property (V) also holds as we only add elements of types $n$ and $n-1$.

In order to verify Property (P) assume that this property fails for a vertex $z'$ of type $n-2$, i.e. there exists a $t$-gon (with $t <m$) formed by vertices of types $n-1$ and $n$ incident with this $z'$. 
By Property (P) for the geometry $(W, \tau, *)$, this $t$-gon has to contain a newly added vertex $x_i$. 
(Note that this vertex has to belong to $\gamma$, as the other newly added vertices are incident with at most one element of type $n$, see item 3 on the list of incidences.)
If $x_i$ is of type $n-1$, then it is incident with at most two vertices of type $n$ (see items 1 and 3 on the list of incidences), one of which is newly added as $m \geq 4$. 
Therefore we may assume that $x_i$ is of type $n$.
The two adjacent vertices of type $n-1$ in $\gamma$ to such a newly added vertex of type $n$ in the $t$-gon are incident with at most one common vertex of type $n-2$, which would be the vertex $z$. 
So we conclude that $z =z'$, for which the construction assures us that Property (P) is still satisfied.

Property (H) is satisfied as we picked a precursor for each newly added vertex of type $n-1$, which we used to define its incidences with elements of types $j\leq n-2$ (see items 2 and 4 on the list of incidences).

Property (C) is satisfied by construction.
\qed

\subsection{Construction of a geometry of type $M$}

Let $\Delta_0 \define (V, \tau, *)$ be a geometry constructed from an infinite countable projective space as in Section~\ref{section:setting}.

Pick an arbitrary vertex $a$ of type $n-1$ in $\Delta_0$. Let $\Lambda_0$ be the geometry obtained from $\Delta_0$ by adding a single vertex of type $n-1$ with precursor $a$, for which the incidences of this added element are completely determined by Property (H). 
One easily verifies that $\Lambda_0$ satisfies the properties listed in Section~\ref{section:inter2}. 

We now define a sequence of geometries $\Lambda_0, \Lambda_1, \dots$ recursively.

Let $S_0$ be an ordered list (indexed by the natural numbers) consisting of those triples $(z,x,y)$ of vertices in $\Lambda_0$ to which one can apply the extension procedure outlined in Section~\ref{section:ext2}. Note that this is possible as the number of such triples is countable.

In step 1 the geometry $\Lambda_1$ is obtained by applying the extension procedure to the first such triple in $S_0$. 
As a consequence additional triples to which one can apply the extension procedure appear, which we order in an ordered list $S_1$ (again indexed by $\mathbb{N}$).

We now proceed by repeating the extension procedure, applied in step $j$ to the first unhandled triple in the list $S_{\nu_2(j)}$, where $\nu_j$ is the 2-adic valuation of $j$, yielding a geometry $\Lambda_j$ and an ordered list of triples $S_j$.

The way we choose wich list to pick a triple from ensures that eventually the extension procedure is applied to every triple appearing on these lists.

Let $\Lambda$ be the direct limit of the sequence $\Lambda_0, \Lambda_1, \dots$. 
We now claim that $\Lambda$ is a geometry of type $M$. 
As this direct limit will still satisfy Properties (F), (H) and (C), one sees that the residues of vertices of type $n-1$ and type $n$ are as intended. 
In particular the residue of a vertex of type $n$ will be isomorphic with the flag geometry of the projective space associated to the vector space $\Gamma$.
The extension procedure applied to every triple together with Property (P), ensures that also the residue of a vertex of type $n-2$ is as intended. 
Residual connectedness is also easily verified as it follows from the previously determined residues of vertices of type $n$, and Property (F). Hence the claim.

\section{Homogeneous geometries of spherical type}\label{section:hom}

In this section we briefly discuss how a combination of the construction in Section~\ref{section:A3} together with the concept of Fra\"iss\'e limits leads to homogeneous geometries. 
The strategy used in this section borrows heavily from~\cite{Ten:11}.

We will restrict our discussion to geometries of types $\mathsf{C}_3$, $\mathsf{H}_3$ and $\mathsf{F}_4$. This allows us to simplify the construction in Section~\ref{section:constr} as we only need to apply Procedure B, and omit Procedures A and C. 

\subsection{Geometries as first order structures}

To start we lay out a first order language $L$ for these geometries and the intermediary objects. 

We have three predicate symbols $V_i$, one for each type $i$ of the Coxeter system $M$, corresponding to the vertices of type $i$ of the geometry. 
We also introduce ternary functions $f_k$ ($k \in \mathbb{N}$), where $f_k(x,y,z)$ is the vertex $a_k$ of the unique shortest path $(a_0 := y, a_1, \dots , a_l := z)$ between $y$ and $z$ in the residue of $x$. 
If such a path does not exist, is not unique or does not make sense (for example if $y$ and $z$ are not contained in the residue of $x$), we set $f_k(x,y,z)$ to be $x$. 
Finally we add binary functions $g_{i,j}$ for each pair of types $i$ and $j$ in $M$ such that $m_{i,j} \geq 4$.
If $x$ and $y$ are two vertices of type $i$ both incident with a vertex $z$ of type $j$, then $g_{i,j}(x,y)$ equals $z$, in each other case $g_{i,j}(x,y)$ equals $x$. 
Note that such an element $z$ is unique by Property D from Section~\ref{section:inter}.

One easily verifies that it is possible to express the axioms of a geometry of type $M$ in this language. 
More exactly we restrict to those geometries satisfying Property (F) from Section~\ref{section:inter}.
By letting the axioms include that every flag is contained in a maximal flag containing a vertex of every type, one automatically has that the geometry is residually connected.

Let $\cC$ be the class of all finitely generated $L$-substructures of those geometries obtained from a finite intermediate geometry in Section~\ref{section:constr} (where we only make use of Procedure B, not of Procedures A and C).

Some examples of $L$-substructures in $\cC$ are flags of arbitrary type (including singular vertices), the Neumaier geometry, etc.

\subsection{The Fra\"iss\'e limit of $\cC$}
We now want to consider the Fra\"iss\'e limit of $\cC$. In order to be able to have this limit we need the following three properties (for details see~\cite[Appendix A]{Fra:00}).
\begin{itemize}
\item 
Hereditary property: If $A \in \cC$, then any finitely generated $L$-substructure of $A$ is again in $\cC$. 
This is clearly satisfied by construction.
\item
Joint embedding property: If $A, B \in \cC$, then there is an $L$-structure $C \in \cC$ such that both $A$ and $B$ are embeddable in $C$. 
\item
Amalgamation property: if $A, B, C \in \cC$ with embeddings $\iota:A \to B$ and $\kappa:A \to C$, then there is a $D \in \cC$ with embeddings $\lambda: B \to D$ and $\mu: C \to D$ such that $\lambda \circ \iota$ equals $\mu \circ \kappa$.
\end{itemize}

In our case the amalgamation property implies the joint embedding property as $\cC$ contains the empty structure. 
The amalgamation property is easily verified by taking the free amalgam of $B$ and $C$ over $A$ (which satisfies the properties listed in Section~\ref{section:inter} by our definition of the language $L$), and then applying the construction of~\ref{section:constr} (again only making use of Procedure B) to obtain a geometry of type $M$ contained in $\cC$.

As all of the properties are satisfied, we are able to consider the Fra\"iss\'e limit $F(\cC)$. 
This will be a countable $L$-structure such that the class of finitely generated $L$-structures embeddable in $F(\cC)$ is exactly $\cC$, and such that each isomorphism between finitely generated $L$-substructures of $F(\cC)$ (which are precisely the structures in $\cC$) can be extended to an automorphism of the entire structure $F(\cC)$.

From the way we set up the language $L$, it is implied that the geometry expressed by the $L$-structure $F(\cC)$ is a geometry of type $M$ which has a automorphism group transitive on flags of the same type (in particular on vertices of the same type), on embedded Neumaier geometries, etc.

\end{document}